\documentclass[11pt,twoside]{article}
\usepackage{cmap} % for russuan cope paste
\def\wdeb{0}
\usepackage{amsfonts}
\usepackage{bezier}
\usepackage{euscript}
\usepackage{amsfonts,amssymb,amsmath,amsbsy,amsthm,amscd}%I included these styles
\usepackage[T2A]{fontenc}
\usepackage[utf8]{inputenc}
\usepackage[english,russian]{babel}   % load Babel setup for English
\usepackage{allerree-utf}
\begin{document}
%======================================================================
\cleanbegin
%======================================================================
\def\udk{515.122.4}
{\LARGE{Краткие сообщения}}
\vskip-0.6cm
\ltitle{Непрерывность обратного в группах}
	{ 
Е.\,А.~Резниченко\footnote[1]{
{\it Резниченко Евгений Александрович} --- канд. физ.-мат. наук, доцент
каф. общей топологии и геометрии мех.-мат. ф-та МГУ,
e-mail: erezn@inbox.ru.\\
\indent {\it Reznichenko Evgenii Aleksandrovich} --- Candidate of Physical and Mathematical Sciences, Associate Professor, Lomonosov
Moscow State University, Faculty of Mechanics and Mathematics, Chair of General Topology and Geometry.
}
}

\iabstract{%АННТОНЦИЯ НА РУССКОМ ЯЗЫКЕ:
Определен широкий класс пространств --- $\Delta$-бэровские пространства.
Паратопологические группы из этого класса являются топологическими.
Класс $\Delta$-бэровских пространств включает локально псевдокомпактные
пространства, бэровские $p$-пространства, бэровские $\Sigma$-пространства и произведения полных по Чеху пространств.
}
{%КЛЮЧЕВЫЕ СЛОВА НА РУССКОМ ЯЗЫКЕ:
паратопологические группы, семитопологические группы,
непрерывность обратного, бэровские пространства, $\Delta$-бэровские пространства.
}
{%АННТОНЦИЯ НА АНГЛИЙСКОМ ЯЗЫКЕ:
We define $\Delta$-Baire spaces. If a paratopological group $G$ is $\Delta$-Baire space,
then $G$ is a topological group. Locally pseudocompact spaces, Baire $p$-spaces, Baire $\Sigma$-spaces, products of \v Cech-complete spaces are $\Delta$-Baire spaces.
}
{%КЛЮЧЕВЫЕ СЛОВА НА АНГЛИЙСКОМ ЯЗЫКЕ:
paratopological groups, semi-topological groups,
continuity of the inverse, Baire spaces, $\Delta$-Baire spaces.}

%------------------------------------------------------------------

Пусть $G$ есть группа и $\tau$ --- топология на $G$.
Группа с топологией $(G,\tau)$ называется {\it семитопологической},
если правые и левые сдвиги непрерывны, иными словами, если умножение
раздельно непрерывно. Группа $(G,\tau)$ называется {\it паратопологической},
если умножение непрерывно. Напомним, группа $(G,\tau)$ называется
{\it топологической}, если умножение и взятие обратного элемента
$g\mapsto g^{-1}$ непрерывны, т.е.\ если $G$ --- паратопологическая группа
с непрерывной операцией взятия обратного элемента.
Начиная с работы Монтгомери [1] 1936 г.\ и на протяжении последующих лет в ряде работ [2--8] рассматривается проблема, которую можно
сформулировать следующим образом: {\it при каких условиях на топологию
паратопологическая группа является топологической группой?} То есть
когда в паратопологической группе операция взятия обратного элемента
непрерывна? В настоящей работе получены дальнейшие продвижения в решении
этой проблемы.
Под пространством подразумеваем регулярное пространство.

Пусть $X$ --- множество, $P\subset X\times X$. Обозначим
$P_x=\{y\in X: (x,y)\in P\}$ для $x\in X$. Если $X$ --- пространство, то $P$
назовем {\it полуоткрытым}, если каждое множество $P_x$ открыто, и $P$ --- {\it полуокрестность диагонали},
если $P$ полуоткрыто и содержит диагональ $\Delta=\{(x,x): x\in X\}$ квадрата $X\times X$.

{\bf Определение.}
Пространство $X$ назовем {\it $\Delta$-бэровским}, если для любой
полуокрестности диагонали $P\subset X\times X$ существует
открытое непустое множество $W\subset X$, для которого $W\times W\subset \overline P$.

Напомним, пространство $X$ называется {\it бэровским}, если пересечение счетного числа плотных
открытых в $X$ множеств непусто и плотно в $X$.
Для группы $G$ и $M\subset G$ обозначим $I(M)=\{(g,h)\in G\times G: g^{-1}h\in M\}$.
Очевидно, $I(M)_g=gM$ для $g\in G$. Поэтому если $G$ --- семитопологическая группа,
а $V$ --- открытая окрестность единицы, то $I(V)$ --- полуокрестность
диагонали.

\begin{theorem}{Теорема 1.}
$\Delta$-Бэровская паратопологическая группа $G$ является топологической группой.
\end{theorem}

{\bf Доказательство.}
Пусть $U$ есть окрестность единицы группы $G$. Для непрерывности
операции взятия обратного в $G$ достаточно найти окрестность $P$ единицы,
для которой $P\subset U^{-1}$.
Существует окрестность $V$ единицы группы $G$, для которой $\overline{V^2}\subset U$.
Покажем, что $\overline{I(V)}\subset I(U)$. Пусть $(g,h)\notin I(U)$, т.е.\ $g^{-1}h\notin U$.
Пусть $O$ есть такая окрестность единицы, что $g^{-1}hO\cap {V^2}=\emptyset$. Тогда
$hO\cap gVV=\emptyset$ и, следовательно, $gV\times hO\cap I(V)=\emptyset$. Поэтому
$(g,h)\notin \overline{I(V)}$. Поскольку $(g,h)\notin I(U)$ мы брали произвольно,
то $\overline{I(V)}\subset I(U)$.
Так как $G$ --- $\Delta$-бэровское пространство, то существует открытое непустое множество
$W\subset G$, для которого $W\times W\subset \overline{I(V)}\subset I(U)$.
Следовательно, множество $P=W^{-1}W\subset U$ и, очевидно, является открытой окрестностью единицы.
Так как $P^{-1}=P$, то $P\subset U^{-1}$.
Теорема доказана.

Прямая Зонгефрея --- бэровская паратопологическая не топологическая группа, не
$\Delta$-бэровское пространство.
Далее наше исследование носит чисто топологический характер,
мы выясняем, какие пространства являются $\Delta$-бэровскими.
Наглядно и просто доказывается, что метризуемое бэровское пространство $X$ является $\Delta$-бэровским.
Приведем схему доказательства. Пусть $P$ --- полуокрестность диагонали $X$.
Для $\varepsilon>0$ положим $M_\varepsilon=\{x\in X: B_\varepsilon(x)\subset P_x\}$, где $B_\varepsilon(x)$ --- $\varepsilon$-окрестность $x$.
Из бэровости $X$ вытекает, что множество $M_\varepsilon$ где-то плотно в $X$ для некоторого $\varepsilon>0$.
Существует открытое непустое множество $W\subset X$, для которого диаметр
$W$ не превосходит $\varepsilon/2$ и множество $W\cap M_\varepsilon$ плотно в $M_\varepsilon$. Тогда
$W\times W\subset \overline P$.
Для обобщения этого факта будем использовать вариации топологической игры
Банаха--Мазура [9,10].
Пусть $X$ есть топологическое пространство, на нем играют два участника
$\alpha$ и $\beta$. Определение топологической игры можно разбить на две части.
Сначала происходит игра в соответствии с некими {\it правилами игры}.
После игры в соответствии с некими {\it правилами определения победителя}
выявляется победитель в игре. Если заданы некие правила игры $T$ и
правила определения победителя $W$, то соответствующую игру будем обозначать
через $G(T,W)$.

{\bf Правила игра $OD$.} Обозначим $U_{-1}=X$.
На $n$-м ходу игрок $\beta$ выбирает открытые непустые множества $V_n,W_n\subset U_{n-1}$,
игрок $\alpha$ выбирает открытое непустое множество $U_n\subset V_n$.

{\bf Правила определения победителя в игре $OD$.}
Игрок $\alpha$ победил, если
$(b)$  последовательность множеств $(W_n)_n$ не локально конечное семейство, т.е.\ она накапливается к некоторой точке;
$(k)$  существует компактное множество $K\subset X$, так что каждое множество $W_n$ пересекает $K$.

Пусть $G$ --- игра на топологическом пространстве $X$, $p$ --- игрок $\alpha$ или $\beta$.
Пространство $X$ называется {\it $p$-благоприятным}
для игры $G$, если у игрока $p$ есть выигрышная стратегия в игре $G$;
$X$ называется {\it $p$-неблагоприятным},
если $X$ не является $p$-благоприятным для игры $G$,
т.е.\ для игрока $p$ не существует выигрышной стратегии в игре $G$.
Стратегия --- это правило (отображение из множества частичных игр), которое
по информации о предыдущих ходах обоих игроков определяет следующий ход игрока.
Также будем называть пространство $X$ $p$-$G$-(не)благоприятным, если оно $p$-(не)благоприятно
для игры $G$.

\begin{theorem}{Теорема 2.}
Если пространство $X$ $\beta$-неблагоприятно для игры $G(OD,b)$,
то $X$ --- $\Delta$-бэровское пространство.
\end{theorem}

{\bf Доказательство.} Предположим противное, т.е.\ существует полуокрестность диагонали $P$, такая, что
$W\times W \setminus \overline P  \ne \emptyset$ для любого непустого открытого множества $W\subset X$.
Определим стратегию для $\beta$ в игре $G(OD,b)$.
На $n$-м шаге игрок $\beta$ выбирает открытые непустые множества
$V_n,W_n\subset U_n$ таким образом, что $V_n\times W_n \cap P=\emptyset$
и $\overline{V_n}\subset U_n$.
Проверим, что определенная стратегия выигрышная для $\beta$, т.е.\ для
получившейся игры $(V_n,W_n,U_n)_n$ не выполняется правило
выигрыша $(b)$ для игрока $\alpha$, а именно семейство $(W_n)_n$ локально конечно.
Пусть $x\in X$. Если $x\notin \bigcap_n V_n$, то $x\notin V_k$ для некоторого $k$
и окрестность $X\setminus \overline{V_{k+1}}$ точки $x$ пересекается с конечным
количеством элементов $(W_n)_n$. Если же $x\in \bigcap_n V_n$,
то окрестность $P_x$ точки $x$ не пересекается ни с одним множеством $W_n$, так
как $x\in V_n$ и $V_n\times W_n \cap P=\emptyset$.
Мы доказали, что $X$ $\beta$-благоприятно для $G(OD,b)$ --- противоречие.
Теорема доказана.

В псевдокомпатном пространстве любая стратегия для игрока $\alpha$ выигрышная в игре $G(OD,b)$,
так что верно

\begin{theorem}{Предложение 1.}
Любое псевдокомпатное пространство $\alpha$-$G(OD,b)$-благоприятно.
\end{theorem}

Пусть $Y\subset X$. Назовем множество $Y$ {\it $C$-плотным}, если для любого счетного семейства $\gamma$ открытых подмножеств $X$ семейство $\gamma$ локально конечно тогда и только тогда, когда семейство $\{U\cap Y:U\in\gamma\}$ локально конечно в $Y$. Для тихоновского пространства $X$ множество $Y$ $C$-плотно в $X$, если и только если множество $Y$ плотно в $X$ и $C$-вложено в $X$.

\begin{theorem}{Предложение 2.}
Если пространство $X$ $\beta$-$G(OD,b)$-неблагоприятно, $Y\subset X$,
(а)  $Y$ плотное $G_\delta$-подмножество $X$,
или (б)  $Y$ --- $C$-плотно в $X$,
или (в)  $Y$ открыто в $X$,
то $Y$ $\beta$-$G(OD,b)$-неблагоприятно.
\end{theorem}

{\bf Доказательство.}
Для открытого в $Y$ множества $U\subset Y$ обозначим $E(U)=X\setminus \overline{Y\setminus U}$. Для случая (а) пусть $Y=\bigcap_n G_n$, где $G_n$ --- открытые плотные подмножества $X$.  Предположим противное, т.е.\ множество $Y$ $\beta$-благоприятно для игры $G(OD,b)$. Пусть $t$ --- выигрышная стратегия для игрока $\beta$. Для случая (в) стратегия $t$ также является выигрышной стратегией для $\beta$ в $G(OD,b)$ на $X$, так что в дальнейшем будем рассматривать случаи (а) и (б).
Определим стратегию $s$ для игрока $\beta$ в игре $G(OD,b)$ на пространстве $X$.
Построим по индукции последовательности открытых множеств $(V_n,W_n,U_n)_n$ пространства $X$ и $(V_n',W_n',U_n')_n$ пространства $Y$. На $n$-м ходу игрок $\beta$ выбирает $U_{n-1}'\subset Y$ таким образом, что
(а) $\overline{U_{n-1}'}\subset U_{n-1}\cap G_n$;
(б) $\overline{U_{n-1}'}\subset U_{n-1}$.
Положим $(V_n',W_n')=t(V_0',W_0',U_0',\dots,V_{n-1}',W_{n-1}',U_{n-1}')$, 
$s(V_0,W_0,U_0,\dots,V_{n-1},W_{n-1},U_{n-1})=(V_n,W_n)=(E(V_n'),E(W_n'))$.
Игрок $\alpha$ выбирает $U_n\subset V_n$.
Игра 
\break
$(V_0',W_0',U_0',\dots)$ выигрышная для $\beta$, поэтому множество $P$ предельных в $X$ точек семейства $(W_n')_n$ не пересекается с $Y$. В случае (а) по построению $P\subset Y$, поэтому множество $P$ пусто. В случае (b)   множество $P$ пусто, поскольку $Y$ $C$-плотно в $X$. Так как $W_n'=W_n\cap Y$ для каждого $n$, то семейство $(W_n)_n$ локально конечно. Следовательно, $s$ --- выигрышная стратегия для $\beta$ и пространство $X$ $\beta$-благоприятно. Противоречие.
Предложение доказано.

Пусть $\mathcal F$ --- семейство подмножеств $X$.
Назовем $\mathcal F$  {\it ограниченным}, если для любого
локально конечного семейства $\gamma$ существует множество $F\in\mathcal F$, для которого 
семейство $\{M\in\gamma: M\cap F\ne\emptyset\}$ конечно.
Для $K\subset X$ будем говорить,
что $\mathcal F$ {\it сходится к} $K$, если для любой окрестности $U\supset K$ семейство 
$\{F\in\mathcal F: F\not\subset U\}$ конечно.

{\bf Правила игра $BM$.} Положим $U_{-1}=X$.
На $n$-м ходу игрок $\beta$ выбирает открытое непустое множество $V_n\subset U_{n-1}$, игрок
$\alpha$ выбирает открытое непустое множество $U_n\subset V_n$.

{\bf Правила определения победителя в игре $BM$.}
Игрок $\alpha$ победил, если
$(i)$  $\bigcap_n U_n\ne\emptyset$;
$(b)$  $(U_n)_n$ --- ограниченное семейство;
$(k)$  $(U_n)_n$ сходится к некоторому компакту $K\subset X$.

Игра $G(BM,i)$ --- это игра Банаха--Мазура. С помощью игры Банаха--Мазура
в [11] получена полезная характеристика бэровских пространств:
пространство $X$ является бэровским, если и только если пространство $X$ $\beta$-неблагоприятно
для игры $G(BM,i)$.
Если $w$ --- одно из определенных
правил выигрыша для игры $OD$ или $BM$, то через $w^*$ будем обозначать
следующее правило выигрыша игрока $\alpha$: либо $\bigcap_n U_n = \emptyset$,
либо выполняется $w$. Если пространство не бэровское, то оно
$\alpha$-благоприятное для $G(R,w^*)$, где $R$ --- либо $OD$, либо $BM$.

\begin{theorem}{Предложение 3.}
Если бэровское пространство $X$ $\alpha$-благоприятно для
игры $G(OD,b^*)$, то $X$ $\beta$-неблагоприятно для $G(OD,b)$.
\end{theorem}

{\bf Доказательство.} В силу теоремы Сан-Ремона [11] достаточно показать, что если пространство 
$X$ $\beta$-благоприятно для $G(OD,b)$ и $\alpha$-благоприятно для
$G(OD,b^*)$, то оно $\beta$-благоприятно для $G(BM,i)$. Пусть --- $t_1$ выигрышная
стратегия для $\beta$ в игре $G(OD,b)$ и $t_2$ --- выигрышная
стратегия для $\alpha$ в игре $G(OD,b^*)$.
Определим выигрышную стратегию $t$ для $\beta$ в игре $G(BM,i)$.

{\bf Первый шаг.}
Пусть
    $(V_0',W_0)=t_1(\emptyset)$,
    $V_0=t_2(V_0',W_0)$.
Положим
    $t(\emptyset)=V_0$.
Пусть $\alpha$ выбирает $U_0\subset V_0$.

{\bf $n$-й шаг.}
Пусть
    $(V_n',W_n)=t_1(V_0',W_0,U_0,\dots,V_{n-1}',W_{n-1},U_{n-1})$,
    $V_n=t_2(V_0',W_0,V_0,\dots,V_{n}',W_{n})$.
Положим
    $t(V_0,U_0,\dots,V_{n-1},U_{n-1})=
        V_n$.
Пусть $\alpha$ выбирает $U_n\subset V_n$.

Проверим, что $t$ --- выигрышная стратегия для $\beta$ в игре $G(BM,i)$,
т.е.\ $\bigcap_n U_n = \emptyset$. В экземпляре
$(V_0',W_0,V_0,\dots,V_n',W_n,V_n,\dots)$ игры $OD$
выиграл $\beta$ по правилам $b$, следовательно,
$\beta$ выиграл по правилам $b$ и в игре
$\xi=(V_0',W_0,U_0,\dots,V_n',W_n,U_n,\dots)$.
С другой стороны, $\alpha$ выиграл в $\xi$ по правилам $b^*$.
Из этого вытекает, что $\bigcap_n U_n = \emptyset$.
Предложение доказано.

\begin{theorem}{Предложение 4.}
Пусть $X$ --- пространство,
$G_1$ и $G_2$ --- одна из следующих пар игр:
$G(BM,k^*)$ и $G(BM,b^*)$;
$G(OD,k^*)$ и $G(OD,b^*)$;
$G(BM,k^*)$ и $G(OD,k^*)$;
$G(BM,b^*)$ и $G(OD,b^*)$.
Тогда если
пространство $X$ $\alpha$-благоприятно для $G_1$, то
оно $\alpha$-благоприятно для $G_2$.
\end{theorem}

{\bf Доказательство.}
Для первых двух пар выигрышная стратегия для $\alpha$ в игре $G_1$ является
выигрышной стратегией для $\alpha$ в $G_2$. В следующих двух парах
игрок $\alpha$ в $G_2$ использует стратегию $\alpha$ в $G_1$, игнорируя $W_n$.
Предложение доказано.

\begin{theorem}{Лемма.}
Пусть $X$ --- пространство,
$(\gamma_n)_n$ --- последовательность семейств открытых множеств,
множество $\bigcup \gamma_n$ плотно в $X$ для каждого $n$. Если для каждой
последовательности $(W_n)_n$, $W_n\in \gamma_n$ с непустым пересечением
а)  $(W_n)_n$ --- ограниченное семейство, то $X$ $\alpha$-благоприятно для игры $G(BM,b^*)$;
б)  $(W_n)_n$ сходится к некоторому компакту, то $X$ $\alpha$-благоприятно для игры $G(BM,k^*)$.
\end{theorem}

{\bf Доказательство.}
Определим для $\alpha$ выигрышную стратегию. На $n$-м шаге будем
выбирать $U_n$ таким образом, что $U_n\subset V_n\cap W_n$ для некоторого
$W_n\in \gamma_n$.
Лемма доказана.

Пространство $X$ называется {\it $p$-пространством}, если существует последовательность $(\gamma_n)_n$  открытых покрытий $X$, таких, что если последовательность $(W_n)_n$, $W_n\in\gamma_n$ имеет непустое пересечение $\bigcap_n W_n$, то $(W_n)_n$ сходится к некоторому компакту. Из леммы, п.б вытекает следующее предложение.

\begin{theorem}{Предложение 5.}
Любое $p$-пространство $\alpha$-благоприятно для игры $G(BM,k^*)$.
\end{theorem}

\begin{theorem}{Предложение 6.}
$\Sigma$-пространство $\alpha$-благоприятно для игры $G(BM,b^*)$.
Сильное $\Sigma$-пространство $\alpha$-благоприятно для игры $G(BM,k^*)$.
\end{theorem}

{\bf Доказательство.}
Пусть $X$ --- (сильное) $\Sigma$-пространство. Тогда существуют $\sigma$-локально конечное семейство $\mathcal F$
замкнутых множеств и покрытие $\mathcal K$, состоящие из замкнутых счетно-компактных (компактных) множеств, так что для каждого $K\in \mathcal K$ и окрестности $U\supset K$ существует множество $F\in \mathcal F$, для которого $K\subset F \subset U$. Пусть $\mathcal F=\bigcup_n \mathcal F_n$, где $\mathcal F_n$ локально конечно и $\mathcal F_n\subset \mathcal F_{n+1}$. Пусть $\gamma_n$ есть семейство открытых $U\subset X$, таких, что для каждого множества $F\in \mathcal F_n$ если $U\cap F\ne \emptyset$, то $U\subset F$. Так как семейство $\mathcal F_n$ локально конечно и состоит из замкнутых множеств, то $\bigcup \gamma_n$ плотно в $X$. Тогда любая последовательность
$(W_n)_n$, $W_n\in\gamma_n$, имеющая непустое пересечение $G=\bigcap_n W_n$, сходится к любому $K\in \mathcal K$, пересекающему $G$.
Если $K$ счетно-компактно, то $(W_n)_n$ ограничено, а если $K$ компактно, то $(W_n)_n$ сходится к компакту. Из леммы  вытекает утверждение предложения.
Предложение доказано.

\begin{theorem}{Предложение 7.}
Пусть $\{X_a: a\in A\}$ --- семейство $\alpha$-$G(OD,k^*)$-благоприятных пространств. Тогда произведение $X=\prod_{a\in A}X_a$ $\alpha$-$G(OD,k^*)$-благоприятно.
\end{theorem}

{\bf Доказательство.} Пусть $t_a$ --- выигрышная стратегия для игрока $\alpha$ в игре $G(OD,k^*)$ на пространстве $X_a$. Построим выигрышную стратегию $t$ для $\alpha$ на $X$. По индукции будем строить конечные попарно непересекающиеся множества $A_n\subset A$ и открытые непустые множества $V_n,W_n,U_n\subset X$,  $V_n^a,W_n^a,U_n^a\subset X_a$ для $a\in A_n$ таким образом, что
(1)
$V_{n+1},W_{n+1}\subset U_n\subset V_n$;
(2)
$U_n=P_n\times \prod_{i=0,\dots,n} \prod_{a\in A_i} U_{n-i}^a,$
где $P_n=\prod_{a\in A\setminus A_n^*}$ и $A_n^*=\bigcup_{i=0,\dots,n} A_i$;
(3)
$V_n^*\subset V_n$ и $W_n^*\subset W_n$, где
$V_n^*=P_n\times \prod_{i=0,\dots,n} \prod_{a\in A_i} V_{n-i}^a$ и
$W_n^*=P_n\times \prod_{i=0,\dots,n} \prod_{a\in A_i} W_{n-i}^a$;
(4)
$V_{n+1}^a,W_{n+1}^a\subset U_n^a\subset V_n^a$ и $U_n^a=t_a(V_0^a,W_0^a,U_0^a,\dots,V_n^a,W_n^a)$ для $a\in A^*=\bigcup_n A_n$.
Положим $A_{-1}=\emptyset$ и $U_{-1}=X$. На $n$-м шаге игрок $\beta$ выбирает открытые непустые множества $V_n,W_n\subset U_{n-1}$. Существуют конечное множество $A_n\subset A\setminus A_{n-1}^*$ и открытые непустые $V_{n-i}^a,W_{n-i}^a\subset X_a$ для $a\in A_i$ и $i=0,\dots,n$, так что выполняется (3). Положим $U_{n-i}^a=t_a(V_0^a,W_0^a,U_0^a,\dots,V_{n-i}^a,W_{n-i}^a)$ для $a\in A_i$ и $i=0,\dots,n$ и определим $U_n$ в соответствии с (2). Положим $t(V_0,W_0,U_0,\dots,V_n,W_n)=U_n$.
Если $\bigcap U_n^a=\emptyset$ для некоторого $a\in A^*$, то $\bigcap U_n=\emptyset$ и игрок $\alpha$ победил. В противном случае из (4) вытекает, что для $a\in A^*$ существует компакт $K^a\subset X_a$, так что $K^a\cap W_n^a\ne\emptyset$ для всех $n$. Пусть $p\in \prod_{a\in A\setminus A^*}$ и $K=\{p\}\times \prod_{a\in A^*}K^a$. Тогда $K\cap W_n\ne\emptyset$ для всех $n$. Мы доказали, что стратегия $t$ выигрышная и пространство $X$ $\alpha$-благоприятно.
Предложение доказано.

\begin{theorem}{Предложение 8.}
Пусть $X$ --- 
$\alpha$-$G(OD,k^*)$-благоприятное пространство и $Y$ --- 
% \hfill \ \ \ \break
%
$\alpha$-$G(OD,b^*)$-бла\-го\-при\-ятное пространство. Тогда произведение $X\times Y$ $\alpha$-$G(OD,b^*)$-благоприятно.
\end{theorem}

{\bf Доказательство.} Пусть $p$ --- выигрышная стратегия для $\alpha$ в $G(OD,k^*)$ на $X$ и $q$ --- выигрышная стратегия для $\alpha$ в $G(OD,b^*)$ на $Y$. Определим стратегию $t$ для $\alpha$ в игре $G(OD,b^*)$ на $X\times Y$. Положим $U^X_{-1}=X$ и $U^Y_{-1}=Y$. На $n$-м шаге игрок $\beta$ выбирает открытые непустые множества $V_n,W_n\subset U_{n-1}= U_{n-1}^X\times U_{n-1}^Y$. Существуют отрытые непустые множества $V_n^X,W_n^X\subset U_{n-1}^X$ и $V_n^Y,W_n^Y\subset U_{n-1}^Y$, такие, что $V_n^X\times V_n^Y\subset V_n$ и $W_n^X\times W_n^Y\subset W_n$. Положим $U_n^X=p(V_0^X,\dots,V_n^X)$, $U_n^Y=q(V_0^Y,\dots,V_n^Y)$ и $t(V_0,\dots,V_n)=U_n=U_n^X\times U_n^Y$.
Проверим, что $t$ --- выигрышная стратегия. Если $\bigcap_n U_n^X$ или $\bigcap_n U_n^Y$ пусто, то $\alpha$ выиграл. В другом случае существует компакт $K\subset X$, который пересекает все $W_n^X$, и последовательность множеств $(W_n^Y)_n$ накапливается к некоторой точке $y\in Y$. Тогда $(W_n)_n$ накапливается к $(x,y)$ для некоторого $x\in K$. 
Предположим, что такой точки $x\in K$ нет. Тогда существует окрестность $O$ точки $y$, для которой $K\times O$ пересекается с конечным числом элементов $(W_n)_n$, и $O$ пересекается с конечным числом элементов $(W_n^Y)_n$ --- противоречие.
Предложение доказано.

Аналогично предложению 2 доказываются следующие два предложения.

\begin{theorem}{Предложение 9.}
Если пространство $X$ $\alpha$-$G(OD,b^*)$-благоприятно, $Y\subset X$,
$Y$ --- плотное $G_\delta$-под\-мно\-жество $X$, или $Y$ $C$-плотно в $X$, или $Y$ открыто в $X$,
то $Y$ $\alpha$-$G(OD,b^*)$-благоприятно.
\end{theorem}

\begin{theorem}{Предложение 10.}
Если пространство $X$ $\alpha$-$G(OD,k^*)$-благоприятно, $Y\subset X$,
$Y$ плотное $G_\delta$-под\-мно\-жество $X$ или $Y$ открыто в $X$, то $Y$ $\alpha$-$G(OD,k^*)$-благоприятно.
\end{theorem}

Из предложений 1, 3--10 вытекает следующая теорема.

\begin{theorem}{Теорема 3.}
Пусть ${\mathcal P}_b$ --- класс всех $\alpha$-$G(OD,b^*)$-благоприятных пространств и ${\mathcal P}_k$ --- класс всех $\alpha$-$G(OD,k^*)$-благоприятных пространств. Тогда
а)
класс ${\mathcal P}_k$ содержит сильные $\Sigma$-пространства и $p$-пространства;
б)
класс ${\mathcal P}_k$ замкнут относительно произведений, перехода к открытым подпространствам, плотным $G_\delta$-подпространствам;
в) ${\mathcal P}_k\subset {\mathcal P}_b$;
г) класс ${\mathcal P}_b$ содержит $\Sigma$-пространства и псевдокомпактные пространства;
д) класс ${\mathcal P}_b$ замкнут относительно перехода к открытым подпространствам, плотным $G_\delta$-подпространствам, $C$-плотным подпространствам;
e) если $X\in {\mathcal P}_k$ и $Y\in {\mathcal P}_b$, то $X\times Y\in {\mathcal P}_b$.
Если $X\in {\mathcal P}_b$ является бэровским пространством, то $X$ $\beta$-$G(OD,b)$-неблагоприятное $\Delta$-бэровское пространство.
Если $X$ гомеомрфно паратопологической группе $G$, то $G$ является топологической группой.
\end{theorem}

%------------------------------------------------------------------------
\spisoklit
{\small\wrefdef{11}

\wref{1}
{Montgomery D.} 
Continuity in topological groups
//
Bull. Amer. Math. Soc. 1936. {\bf 42.}  879--882.

\wref{2}
{Ellis R.}
A note on the continuity of the inverse
//
Proc. Amer. Math. Soc. 1957. {\bf 8.}  372--373.

\wref{3}
{Zelazko W.}
A theorem on $B_0$ division algebras
//
Bull. Acad. Pol. Sci. 1960. {\bf 8.} 373--375.

\wref{4}
{Brand N.}
Another note on the continuity of the inverse
//
Arch. Math. 1982. {\bf 39.} 241--245.

\wref{5}
{Pfister H.}
Continuity of the inverse
//
Proc. Amer. Math. Soc. 1985. {\bf 95.} 312--314.

\wref{6}
{Reznichenko E.A.}
Extensions of functions defined on products of pseudocompact spaces and continuity of the inverse in pseudocompact groups
//
Topol. Appl. 1994. {\bf 59.}  233--244.

\wref{7}
{Bouziad A.}
Every \v Cech-analytic Baire semitopological group is a topological group
//
Proc. Amer. Math. Soc. 1998. {\bf 24,} N3. 953--959.

\wref{8}
{Arhangel'skii A.V.,  Reznichenko E.A.}
Paratopological and semitopological groups versus topological groups
//
Topol. Appl. 2005. {\bf 151.} 107-119.

\wref{9}
{Choquet G.}
Lectures on Analysis, Vol. I. N. Y.: W.~A.~Benjamin Inc, 1969.

\wref{10}
{Telgarsky R.}
Topological games: On the 50th anniversary of the Banach--Mazur game
//
Rocky Mount. J. Math. 1987. {\bf 17.}  227--276.

\wref{11}
{Saint Raymond J.}
Jeux topologiques et espaces de Namioka
//
Proc. Amer. Math. Soc. 1983. {\bf 87.} 499--504.
}

\bigskip
\hfill
\parbox{11em}{
Поступила в редакцию 
\\
06.09.2021
}
\smallskip
\lend

\end{document}